\input amstex
\input amsppt.sty
\magnification=\magstep1
\hsize=30truecc
\baselineskip=16truept
\vsize=22.2truecm
\TagsOnRight
\nologo
\pageno=1
\topmatter
\def\N{\Bbb N}
\def\Z{\Bbb Z}

\def\l{\left}
\def\r{\right}
\def\bg{\bigg}
\def\({\bg(}
\def\[{\bg[}
\def\){\bg)}
\def\]{\bg]}
\def\t{\text}
\def\f{\frac}
\def\colon{{:}\;}
\def\mo{\roman{mod}}
\def\em{\emptyset}
\def\se {\subseteq}
\def\sp {\supseteq}
\def\sm{\setminus}

\def\bi{\binom}
\def\eq{\equiv}
\def\cs{\cdots}
\def\ls{\leqslant}
\def\gs{\geqslant}
\def\al{\alpha}

\def\h{[G:H]}
\def\hk{H\cap K}

\def\ik{^k_{i=1}}

\def\Proof{\noindent{\it Proof}}
\def\Remark{\medskip\noindent{\it  Remark}}
\def\Def{\medskip\noindent{\it  Definition}}
\def\Ack{\medskip\noindent {\bf Acknowledgment}}
\hbox{Internat. J. Math. 17(2006), no.\,9, 1047--1064.}
\bigskip
\title Finite covers of groups by cosets or subgroups\endtitle
\author Zhi-Wei Sun\endauthor
\affil Department of Mathematics and Institute of Mathematical Science
\\Nanjing University, Nanjing 210093, People's Republic of China
    \\ {\tt zwsun\@nju.edu.cn}
    \\ {\tt http://pweb.nju.edu.cn/zwsun}
\endaffil
\abstract This paper deals with combinatorial aspects
of finite covers of groups by cosets or subgroups. Let $a_1G_1,\ldots,a_kG_k$
be left cosets in a group $G$ such that $\{a_iG_i\}_{i=1}^k$
covers each element of $G$ at least $m$ times but none of its proper
subsystems does. We show that if $G$ is cyclic, or $G$ is finite and
$G_1,\ldots,G_k$ are normal Hall subgroups of $G$, then the inequality
$k\gs m+f([G:\bigcap_{i=1}^kG_i])$ holds, where
$f(\prod_{t=1}^rp_t^{\al_t})=\sum_{t=1}^r\al_t(p_t-1)$ if
$p_1,\ldots,p_r$ are distinct primes
and $\al_1,\ldots,\al_r$ are nonnegative integers.
When all the $a_i$ are the identity element of $G$
and all the $G_i$ are subnormal in $G$,  we prove that
there is a composition series from $\bigcap_{i=1}^kG_i$ to $G$
whose factors are of prime orders. The paper also includes some other results
and two challenging conjectures.
\endabstract

\thanks 2000 {\it Mathematics Subject Classification.} Primary 20D60;
Secondary 05A05, 11B25.
\newline
\indent
The author is supported by the National Science Fund
for Distinguished Young Scholars (No. 10425103) and a Key
Program of NSF in P. R. China.
\endthanks
\endtopmatter
\document

\heading {1. Introduction}\endheading
Let $G$ be a multiplicative group, and let
$$\Cal A=\{a_iG_i\}^k_{i=1}\tag1.1$$
be a finite system of left cosets in $G$
(where $a_1,\ldots,a_k\in G$, and $G_1,\ldots,G_k$ are subgroups of $G$).

For any $I\se[1,k]=\{1,\ldots,k\}$
we define the {\it index map} $I^*$ from $G$
to the power set of $[1,k]$ as follows:
$$I^*(x)=\{i\in I\colon  x\in a_iG_i\}\quad (x\in G).\tag1.2$$
$m(\Cal A)=\inf_{x\in G}|[1,k]^*(x)|$ is said to be
the {\it covering multiplicity} of (1.1).

For a positive integer $m$, if $m(\Cal A)\gs m$ (respectively, $|[1,k]^*(x)|=m$
for all $x\in G$)
then we call (1.1) an {\it $m$-cover} (resp. {\it exact $m$-cover}) of $G$.
If (1.1) forms an $m$-cover of $G$
but none of its proper subsystems does,
then it is said to be a {\it minimal} (or an {\it irredundant})
$m$-{\it cover} of $G$.
A {\it cover} of $G$ refers to a 1-cover of $G$, and
a {\it disjoint cover} (or {\it partition}) of $G$
means an exact 1-cover of $G$. Obviously an exact $m$-cover is
 a minimal $m$-cover and
any minimal $m$-cover has covering multiplicity $m$.
Disjoint covers of general groups were investigated by I. Korec [K74],
M. M. Parmenter [P84] and R. Brandl [Br90],
while exact $m$-covers of groups were studied by the author
in [S01] and [S04b].

If (1.1) forms a minimal cover of a group $G$,
then all the $G_i$ have finite index in $G$
by B. H. Neumann [N54a, N54b], and furthermore
$$\[G:\bigcap_{i=1}^kG_i\]\ls k!\tag1.3$$
by M. J. Tomkinson [T87].
The author [S90] observed that this still holds if
$(1.1)$ is a minimal $m$-cover of a group $G$.

 Here we give a simple way to explain why $[G:G_i]<\infty$
for all $i=1,\ldots,k$ provided that (1.1) is a minimal $m$-cover of a group
$G$. For any $j\in[1,k]$, there is an $x\in a_jG_j$ with $|[1,k]^*(x)|=m$.
Choose a minimal $I\se([1,k]\sm[1,k]^*(x))\cup\{j\}$ such that
$\{a_iG_i\}_{i\in I}$ is a cover of $G$, then we must have
$j\in I$ and hence $[G:G_j]<\infty$ by Neumann's result.

 Neumann's basic result has applications in Galois theory, group rings,
Banach spaces, projective geometry and Riemann surfaces (cf. [SV]).

 Concerning covers of groups by infinitely many cosets, Tomkinson [T86]
 proved that if a group $G$ is irredundantly
covered by $\kappa\gs
\aleph_0$ left cosets $a_iG_i$ then $[G:G_i]\ls2^{\kappa}$
for each $i$ and hence
$$\[G:\bigcap_iG_i\]\ls\prod_{i}[G:G_i]\ls(2^{\kappa})^{\kappa}=2^{\kappa}.$$

Any infinite cyclic group is isomorphic to the additive group $\Z$
of integers. The subgroups of $\Z$ different from $\{0\}$ are in
the form $n\Z=\{nx:\, x\in\Z\}$ where $n\in\Z^+=\{1,2,3,\ldots\}$.
For any positive integer $n$, the index of $n\Z$ in $\Z$ is $n$
and a coset of $n\Z$ in $\Z$ is just a residue class
$$a+n\Z=\{x\in\Z:\, x\eq a\ (\mo\ n)\}\quad (a\in\Z).$$
The study of cover of $\Z$ in the form
$$A=\{a_i+n_i\Z\}^k_{i=1}\tag1.4$$
was initiated by
P. Erd\"os ([E50]) in the early 1930s, who viewed
this topic as his most favorite one (cf. [E97]).
Since then many researchers have investigated such covers
of $\Z$ (cf. [G04] and [PS]). A simple example of cover of $\Z$ with distinct
moduli is
$$\{2\Z,\ 3\Z,\ 1+4\Z,\ 5+6\Z,\ 7+12\Z\}.$$
M. Z. Zhang [Z91] showed that for each $m=2,3,\ldots$
there are infinitely many exact $m$-covers of $\Z$
which cannot be split into an exact $n$-cover of $\Z$ (with $0<n<m$) and
an exact $(m-n)$-cover of $\Z$.
The author investigated
$m$-covers of $\Z$ and exact $m$-covers of $\Z$ in
a series of papers (see, e.g., [S97, S99, S03a, S04a]).
Covering multiplicity plays an important role in the author's unification
of zero-sum problems, subset sums and covers of $\Z$ (cf. [S03b]).

Let $H$ be a subnormal subgroup of a group $G$ with finite index. Define
$$d(G,H)=\sum_{i=1}^n(|H_i/H_{i-1}|-1)\tag1.5$$
where $H=H_0\subset  H_1\subset \cs\subset H_n=G$
is a composition series from $H$ to $G$. ($d(G,G)$ is regarded as $0$.)
By [S90, Theorem 6] we have
$$[G:H]-1\gs d(G,H)\gs f([G:H])\gs \log_2[G:H],\tag1.6$$
here the Mycielski function $f:\Z^+\to\N=\{0,1,2,\ldots\}$
is given by
$$f\(\prod_{t=1}^rp_t^{\al_t}\)=\sum^r_{t=1}\al_t(p_t-1)\tag1.7$$
where $p_1,\ldots,p_r$ are distinct primes and
$\al_1,\ldots,\al_r\in\N$. In [S01] the author showed that
$d(G,H)=f([G:H])$ if and only if $G/H_G$ is solvable,
where $H_G=\bigcap_{g\in G}gHg^{-1}$
is the largest normal subgroup of $G$ contained in $H$.

Here is the main result of [S01]
on exact $m$-covers of groups.

\proclaim{Theorem 1.1 {\rm (Z. W. Sun, 2001)}}
Let $(1.1)$ be an exact $m$-cover of a group
$G$ by left cosets.

{\rm (i)} Whenever $G/(G_i)_G$ is solvable, we have $k\gs
m+f([G:G_i])$.

{\rm (ii)} If all the $G_i$ are subnormal in $G$, then we have the
inequality
$$k\gs m+d\(G,\bigcap_{i=1}^kG_i\),\tag1.8$$
where the lower bound is best possible.
\endproclaim

We now give some historical remarks. When $m=1$ and $G$ is
abelian, part (i) was first conjectured by J. Mycielski (cf. [MS])
and confirmed by \v S. Zn\'am [Zn66] in the case $G=\Z$. In 1988
M. A. Berger, A. Felzenbaum and A. S. Fraenkel [BFF] proved that
if $(1.1)$ forms a disjoint cover of a finite solvable group $G$
then $k\gs 1+f([G:G_i])$ for any $i=1,\ldots,k$. When $m=1$ and
$G_1,\ldots,G_k$ are normal in $G$, part (ii) was essentially
obtained by I. Korec [K74].

For general covers of groups, things become
more subtle and complicated. By using
algebraic number theory and characters of abelian groups,
G. Lettl and the author [LS] recently obtained the following
result.

\proclaim{Theorem 1.2 {\rm (G. Lettl and Z. W. Sun, 2004)}}
Let $(1.1)$ be a minimal $m$-cover of an abelian group $G$.
Then $k\gs m+f([G:G_i])$ for all $i=1,\ldots,k$.
\endproclaim

As Example 1.1 of [LS] shows, we cannot replace the lower bound of
$k$ in Theorem 1.2 by $m+f([G:\bigcap_{i=1}^kG_i])$ even if
$G=C_p\times C_p$ where $p$ is an odd prime and $C_p$ is the
cyclic group of order $p$.

\Def\ 1.1. Let (1.1) be a finite system of left cosets
in a group $G$. If $[1,k]^*(G)=\{[1,k]^*(x)\colon x\in G\}$
contains $I^*(G)$ for no $\em\not=I\subset [1,k]$
(i.e., whenever $\em\not=I\subset[1,k]$ there is a $g\in G$
such that $I^*(g)\not=[1,k]^*(x)$ for all $x\in G$),
then we call (1.1) a {\it regular system}.
If (1.1) is regular, and $[1,k]^*(G)\not\supseteq\em^*(G)=\{\em\}$
(i.e., (1.1) is a cover of $G$), then we call (1.1)
{\it a regular cover} of $G$.
\medskip

When (1.1) is a minimal $m$-cover of a group $G$,
it forms a regular cover of $G$ because
for any $I\subset[1,k]$ there is a $g\in G$
such that $|I^*(g)|<m$ while $|[1,k]^*(x)|\gs m$ for all $x\in G$.
A regular system may not be a cover, e.g., we can
easily check that $\{2\Z,4\Z,2+4\Z\}$ is a regular system of residue classes
but it is not a cover of $\Z$.
As we will see later, if (1.1) forms a regular cover of a group $G$
then all the indices $[G:G_i]$ must be finite.

Now we introduce our first result in this paper.

\proclaim {Theorem 1.3} If $(1.1)$ is a regular cover of a group $G$,
and for any $i,j=1,\ldots,k$ either $G_i$ and $G_j$
are subnormal in $G$ with $[G:G_i]$ relatively prime to
$[G_i:G_i\cap G_j]$, or $G_i$ and $G_j$ are normal in $G$ with $G/(G_i\cap G_j)$
cyclic, then we have the inequality
$$k\gs m(\Cal A)+d\(G,\bigcap^k_{i=1}G_i\).\tag1.9$$
\endproclaim

\Remark\ 1.1. (a) In view of Example 1.1 of [LS], the conditions of Theorem 1.3
are essentially
indispensable. (b) By Example 1.2 of [S01],
for any subnormal subgroup $H$ of $G$
with finite index, there does exist an exact $m$-cover
(1.1) of $G$ such that all the $G_i$ are subnormal in $G$,
$\bigcap_{i=1}^kG_i=H$ and $k=m+d(G,H)$.
\medskip

A subgroup $H$ of a finite group $G$ is called a {\it Hall subgroup} of $G$
if $|H|$ is relatively prime to $[G:H]$ (i.e., $1$ is the greatest
common divisor $(|H|,[G:H])$ of $|H|$ and $[G:H]$). By [S04b, Corollary 2.1],
a subnormal Hall subgroup of $G$ must be normal in $G$.

\proclaim{Corollary 1.1} Let $(1.1)$ be a regular cover of a group $G$.
If $G$ is cyclic, or $G$ is finite and $G_1,\ldots,G_k$
are normal Hall subgroups of $G$, then $(1.9)$ holds.
\endproclaim
\Proof. This follows from Theorem 1.3 immediately. \qed

\Remark\ 1.2. The following consequence of Corollary 1.1
was announced by the author in [S01]:
Let (1.1) be a minimal $m$-cover of a group $G$. If
$G$ is cyclic, or $|G|$ is squarefree and all the $G_i$ are normal in $G$,
then we have the inequality
$$k\gs m+f\(\[G:\bigcap_{i=1}^kG_i\]\).\tag1.10$$
When $m=1$ and $G=\Z$, this was conjectured by Zn\'am [Zn75] and
confirmed by R. J. Simpson [Si85]; when $m=1$ and $G$ is of squarefree order,
this was proved by Berger, Felzenbaum
and Fraenkel [BFF] in 1988.
\medskip

Recall that a group $G$ is said to be {\it perfect} if it coincides with
its derived group $G'$. Here is our second theorem.

\proclaim{Theorem 1.4} Suppose that  $\{G_i\}_{i=1}^k$ is
a minimal $m$-cover of a group $G$ by subnormal subgroups.
Then there is a composition series from $\bigcap_{i=1}^kG_i$ to $G$
whose factors are of prime orders, and all the $G_i$ contain
every perfect subgroup of $G$.
\endproclaim

\Remark\ 1.3. (a) Clearly Theorem 1.4 extends the following result of
M. A. Brodie, R. F. Chamberlain and L.-C. Kappe [BCK]:
If $\{G_i\}_{i=1}^k$ is a minimal cover of a group $G$
by finitely many normal subgroups, then
$G/\bigcap\ik G_i$ is solvable and all perfect normal subgroups of $G$
are contained in each of $G_1,\ldots,G_k$.
(b) Any finite non-cyclic group $G$ can be covered by finitely many
proper subgroups because $G=\bigcup_{x\in G}\langle x\rangle$,
but no field can be covered by finitely many proper subfields
(cf. [BBS]).

\medskip
{\it Example}\ 1.1. Let $G$ be a group with $G/Z(G)$ finite, where
$Z(G)$ is the center of $G$.
For $x,y\in G$, if $xy\ne yx$ then
$(x^{-1}y)x\not=x(x^{-1}y)$ and hence $x Z(G)\ne yZ(G)$.
Let $X=\{x_1,\ldots,x_k\}$ be a maximal
set of pairwise non-commuting elements of $G$. Then $k=|X|\ls|G/Z(G)|$,
and $\{C_G(x_i)\}^k_{i=1}$ forms a minimal cover of $G$
by centralizers with
$\bigcap^k_{i=1}C_G(x_i)=Z(G)$ (cf. [T87, Theorem 5.1]), and
$|G/Z(G)|\ls c^k$ for some absolute constant $c>0$ (see [Py87]).
If $\overline G=G/Z(G)$ is not solvable,
then not all the $C_G(x_i)$ are subnormal in $G$ by Theorem 1.4.
When $\overline G$ is solvable, Tomkinson [T97] provided a lower bound of $k$
(conjectured by Cohn [C94])
in terms of a chief factor of $\overline G$. D. R. Mason [M78] proved that
$|G|\gs2k-2$, which was conjectured by Erd\"os and E. G. Straus [ES].
\medskip

Concerning covers of a group by subnormal subgroups,
the author and his student Song Guo have made the following conjecture.

\proclaim{Conjecture 1.1 {\rm (S. Guo and Z. W. Sun, 2004)}} Let
$\{G_i\}^k_{i=1}$ be a minimal $m$-cover of
a group $G$ by finitely many subnormal subgroups.
Assume that $[G:\bigcap^k_{i=1}G_i]=\prod^r_{t=1}p_t^{\al_t}$,
where $p_1,\ldots,p_r$ are distinct primes and $\al_1,\ldots,\al_r$
are positive integers. Then
$$k>m+\sum^r_{t=1}(\al_t-1)(p_t-1).$$
\endproclaim

Let $H$ be a subgroup of a group $G$. When
$X\se G$ is a union of finitely many left cosets of $H$,
we use $[X:H]$ to represent the number of left cosets of $H$ contained
in $X$, which agrees with the index of $H$ in $X$
if $X$ is a subgroup of $G$. Below is our third result in this paper.

\proclaim{Theorem 1.5} Let $H$ be a subgroup of a group $G$ with
finite index. Suppose that $G_1,\ldots,G_k$ are subnormal subgroups of $G$
containing $H$ with $([G:G_i],[G_i:H])=1$ for all $i=1,\ldots,k$.
Given $a_1,\ldots,a_k\in G$ we have the inequality
$$\[\bigcup_{i=1}^ka_iG_i:H\]\gs\[\bigcup_{i=1}^kG_i:H\].\tag1.11$$
\endproclaim

Here is an obvious consequence of Theorem 1.5.

\proclaim{Corollary 1.2} Let $G$ be a finite group and let
$G_1,\ldots,G_k$ be normal Hall subgroups of $G$. Then, for any
$a_1,\ldots,a_k\in G$, we have
$$\bg|\bigcup_{i=1}^ka_iG_i\bg|\gs\bg|\bigcup_{i=1}^kG_i\bg|.\tag1.12$$
\endproclaim

\Remark\ 1.4. In [S90] it was asked whether for any
finitely many left cosets $a_1G_1,\ldots,a_kG_k$ in a finite group $G$
we always have (1.12). Later Tomkinson pointed out that this is not true
for the Klein group $C_2\times C_2$.
\medskip

 The following conjecture of the author seems very challenging.

\proclaim{Conjecture 1.2 {\rm (Z. W. Sun, 2004)}} Let
$a_1G_1,\ldots,a_kG_k\ (k>1)$ be finitely many pairwise disjoint left cosets in
a group $G$ with $[G:G_i]<\infty$ for all $i=1,\ldots,k$. Then
$([G:G_i],[G:G_j])\gs k$ for some $1\ls i<j\ls k$.
\endproclaim

 \Remark\ 1.5. (a) The case $k=2$ will be handled in Remark 2.2.
 The conjecture remains open even for the
 additive cyclic group $\Z$.

 (b) Under the condition of Conjecture 1.2, clearly
 $$\bg[G:\bigcap_{j=1}^kG_j\bg]
 \gs\bg[\bigcup_{i=1}^ka_iG_i:\bigcap_{j=1}^kG_j\bg]
 =\sum_{i=1}^k\bg[G_i:\bigcap_{j=1}^kG_j\bg]$$
 and hence $\sum_{i=1}^k[G:G_i]^{-1}\ls 1$.
 Suppose that $[G:G_1]\ls\cdots\ls [G:G_k]$.
 Since $\sum_{i=1}^{k-1}[G:G_i]^{-1}<1$, there is an
 $i\in[1,k-1]$ such that $[G:G_i]\not\ls k-1$ and hence $[G:G_k]\gs[G:G_i]\gs k$.
 If $[G:G_k]$ is divisible by all those $[G:G_1],\ldots,[G:G_{k-1}]$
 (this happens if $G$ is a $p$-group with $p$ a prime), then
 $([G:G_i],[G:G_k])=[G:G_i]\gs k$.

 \medskip

 We will derive some combinatorial properties of regular systems
in the next section,
and prove Theorems 1.3-1.5 in the third section.

\heading{2. Combinatorial Properties of Regular Systems}\endheading

  Recall that a {\it monoid} is a semigroup containing an identity element.

  In order to unify disjoint cover and minimal cover, the author [S90] proposed
the following general notion.

Let $M$ be a commutative monoid (considered as an additive one) and $S$ a set.
A finite system $\{X_i\}_{i=1}^k$
of nonempty subsets of a set $X\not=\em$ is called
an $(M,S)$-{\it cover}
of $X$ if there exist $m_1,\ldots,m_k\in M$ such that
$$\bigg\{\sum^k\Sb i=1\\x\in X_i\endSb m_i:\, x\in  X \bigg\}\se S$$
but
$$ \bigg\{\sum \Sb i\in I\\x\in  X_i\endSb
m_i:\, x\in  X\bigg\}\not\se S\ \qquad
\t{for all}\ I\subset [1,k].$$
(Without any loss of generality we may let $S$ be a subset of $M$.)

  Though useful in [S90], this concept involving a monoid and a set
seems cumbersome. We can avoid it by using regular covers
of a nonempty set by nonempty subsets. (The definition is similar to that
of regular cover of a group by cosets.)

Let $X$ be a nonempty set and $\{X_i\}_{i=1}^k$
a finite system of nonempty subsets of $X$.
If $\{X_i\}_{i=1}^k$ forms an $(M,S)$-cover of $X$ for some commutative
monoid $M$ and its subset $S$, then it is easy to see that
$\{X_i\}_{i=1}^k$ is a regular cover of $X$.
Conversely, if $\{X_i\}_{i=1}^k$ is a regular cover of $X$
then it is an $(M_0,S_0)$-cover of $X$,
where $M_0$ is the multiplicative monoid of positive integers,
and $$S_0=\bigg\{\prod_{i\in  [1,k]^*(x)}p_i\colon  x\in  X\bigg\}$$
with $p_1,\ldots,p_k$ distinct primes.

\proclaim{Theorem 2.1} Let $(1.1)$ be a regular system of
left cosets in a group $G$.

{\rm (i)} Let $\em\ne I\se[1,k]$, and $I\not=[1,k]$ if $(1.1)$ is not a cover of $G$.
 Set $\bar I=[1,k]\sm I$, and regard $\bigcap_{i\in\em}G_i$ as $G$.

\ \ {\rm (a)} $\bigcup_{i\in  I} a_iG_i$ contains
a left coset of $\bigcap_{j\in\bar I}G_j$,
and $\bigcap_{j\in\bar I}G_j$ contains
at most $|I|!$ left cosets of $\bigcap^k_{i=1}G_i$.
For each $x\in\bigcap_{j\in\bar I}G_j$
there are positive integers $n_1\ls 1,\ldots,n_{|I|}\ls |I|$
such that $x^{n_1\cs n_{|I|}}\in\bigcap^k_{i=1}G_i$.

\ \  {\rm (b)} If $[G:G_i]<\infty$ for all $i=1,\ldots,k$, then
$$ \sum_{i\in  I}\f1{[G:G_i]}\gs\sum_{i\in I}
\f1{[G:G_i\cap\bigcap_{j\in\bar I}G_j]}\gs \f1{[G:\bigcap_{j\in \bar I}G_j]}
\gs \f1{\prod_{j\in\bar I}[G:G_j]}\tag2.1$$
where an empty product is regarded as $1$. Also,
$$\[\bigcup_{i\in I}G_i:\bigcap_{i=1}^kG_i\]\ls\sum^{|I|}_{l=1}
(-1)^{l-1}(k-l)!\bi{|I|}l.\tag2.2$$

 {\rm (ii)} For any subgroup $H$ of $G$ with
 $J=\{1\ls i\ls k\colon G_i\not\se H\}\not=\em,[1,k]$,
we have
$$|\{C\in G/H:\,C\sp a_iG_i\ \t{for some}\ 1\ls i\ls k\}|
\gs\[\(\bigcap_{j\in J} G_j\)H:H\],
\tag2.3$$
where $G/H$ refers to the set $\{gH:\,g\in G\}$.
\endproclaim
\Proof. (i) Let us prove part (i) first.

(a) Since (1.1) is regular and $\bar I\not=[1,k]$,
there exists a $g\in G$ such that
$${\bar I}^*(g)=\{j\in\bar I:\,g\in a_jG_j\}\not\in[1,k]^*(G).$$
For each $x\in\bigcap_{j\in\bar I}G_j$,
as ${\bar I}^*(gx)={\bar I}^*(g)\ne[1,k]^*(gx)$,
we must have
$gx\in\bigcup_{i\in I}a_iG_i$. So
$$g\(\bigcap_{j\in\bar I}G_j\)\se\bigcup_{i\in I}a_iG_i.$$
For each $i\in I$, $\{i\}\cup\bar I$ is the complement of
 $I\sm\{i\}$ in $[1,k]$;
if $a_iG_i\cap g(\bigcap_{j\in\bar I}G_j)$ is nonempty
then it contains exactly $[G_i\cap\bigcap_{j\in\bar I}G_j:F]$
left cosets of $F=G_1\cap\cdots\cap G_k$.
As
$$g\(\bigcap_{j\in\bar I}G_j\)
=\bigcup_{i\in I}\(a_iG_i\cap g\(\bigcap_{j\in\bar I}G_j\)\),$$
we have
$$\bigg[\bigcap_{j\in\bar I} G_j:F\bigg]
\ls\sum_{i\in I}\bigg[\bigcap_{j\in\overline{I\sm\{i\}}}G_j:F\bigg].\tag2.4$$
In the case $I=\{i\}$, this yields that $[\bigcap_{j\in\bar
I}G_j:F]\ls1=|I|!$.
By induction on $|I|$ we see that $[\bigcap_{j\in\bar I}G_j:F]\ls|I|!$.

Let $x\in\bigcap_{j\in\bar I}G_j$. Since
$$\l\{g,gx,\ldots,gx^{|I|}\r\}\se g\(\bigcap_{j\in\bar I}G_j\)
\se\bigcup_{i\in I}a_iG_i,$$
for some
$0\ls s<t\ls|I|$ both $gx^s$ and $gx^t$ belong to a certain
$a_iG_i$ with $i\in I$
and therefore
$x^{t-s}=(gx^s)^{-1}gx^t\in G_i$. This shows that
$$x^n\in\bigcup_{i\in I}G_i\cap\bigcap_{j\in\bar I}G_j=\bigcup_{i\in I}
\bigcap_{j\in\overline{I\sm\{i\}}}G_j$$
where $n=t-s\in\{1,\ldots,|I|\}$.
Recall that  if $|I|=1$ then $\bigcap_{j\in\bar I}G_j=F$.
Again, by induction on $|I|$, we find that $x^{n_1\cs n_{|I|}}\in F$
for some positive integers $n_1\ls 1,\ldots,n_{|I|}\ls|I|$.

(b) When $G_1,\ldots,G_k$ are of finite index in $G$,
dividing both sides of the inequality $(2.4)$ by $[G:F]<\infty$ we obtain (2.1).

For each nonempty subset $J$ of $I$, if $J=[1,k]$ then
$$\[\bigcap_{j\in J}G_j:F\]-\sum\Sb i\in I\\i>\max J\endSb
\[G_i\cap\bigcap_{j\in J}G_j:F\]=1-0=(k-|J|)!$$
where $\max J$ is the maximal element of $J$;
otherwise, by $(2.4)$ and part (a) we have
$$\aligned&\[\bigcap_{j\in J}G_j:F\]-\sum\Sb i\in I\\i>\max J\endSb
\[G_i\cap\bigcap_{j\in J}G_j:F\]
\\\ls&\sum^k\Sb i=1\\i\not\in J\endSb\[G_i\cap\bigcap_{j\in J}G_j:F\]
-\sum\Sb i\in I\\i>\max J\endSb
\[G_i\cap\bigcap_{j\in J}G_j:F\]
\\=&\sum\Sb1\ls i\ls k,\ i\not\in J\\i\not\in I\ \t{or}\ i\ls\max J\endSb
\[G_i\cap\bigcap_{j\in J}G_j:F\]
\\\ls&\sum\Sb1\ls i\ls k,\ i\not\in J\\i\not\in I\ \t{or}
\ i\ls\max J\endSb(k-(|J|+1))!
\\=&\sum^k\Sb i=1\\i\not\in J\endSb(k-(|J|+1))!
-\sum\Sb i\in I\\i>\max J\endSb(k-(|J|+1))!
\\=&(k-|J|)!-\sum\Sb i\in I\\i>\max J\endSb(k-(|J|+1))!.
\endaligned$$

If $\em\not=J'\se I$ and $2\mid|J'|$, then
$J'=J\cup\{i\}$ and $2\nmid|J|$
where $i=\max J'$ and $J=J'\sm\{i\}$.
Thus, in view of the above and the inclusion--exclusion principle,
$$\aligned\[\bigcup_{i\in I}G_i:F\]
=&\sum_{\em\ne J\se I}(-1)^{|J|-1}\[\bigcap_{j\in J}G_j:F\]
\\=&\sum\Sb J\se I\\2\nmid|J|\endSb\(\[\bigcap_{j\in J}G_j:F\]
-\sum\Sb i\in I\\i>\max J\endSb
\[G_i\cap\bigcap_{j\in J}G_j:F\]\)
\\\ls&\sum\Sb J\se I\\2\nmid|J|\endSb \((k-|J|)!
-\sum\Sb i\in I\\i>\max J\endSb (k-(|J|+1)!\)
\\=&\sum_{\em\not=J\se I}(-1)^{|J|-1}(k-|J|)!
=\sum^{|I|}_{l=1}(-1)^{l-1}\bi{|I|}l(k-l)!.
\endaligned$$
So (2.2) is valid.

(ii) By part (i) there exists $g\in G$ such that
$g(\bigcap_{j\in J}G_j)\se\bigcup_{i\in\bar J}a_iG_i$.
Therefore
$$g\(\bigcap_{j\in J}G_j\)H\se \bigcup_{i\in\bar J}a_iG_iH
=\bigcup_{i\in\bar J}a_iH$$
and hence
$$\align&|\{C\in G/H:\, C\supseteq a_iG_i\ \t{for some}\ 1\ls i\ls k\}|
\\&\qquad=|\{a_iH\colon 1\ls i\ls k\ \t{and}\ a_iG_i\se a_iH\}|
\\&\qquad\gs\[\(\bigcap_{j\in J}G_j\)H:H\].
\endalign$$
We are done. \qed

\Remark\ 2.1. (i) By Theorem 2.1(i), for any regular cover (1.1) of
a group $G$ we have
 $$\[G:\bigcap^k_{i=1}G_i\]=\[\bigcap_{i\in\em}G_i:\bigcap_{i=1}^kG_i\]
 \ls k!<\infty.$$
(ii) Let (1.1) be a minimal cover of a group $G$.
That $\sum_{i\in I}[G:G_i]^{-1}\gs\prod_{j\in\bar I}[G:G_j]^{-1}$
for all $\em\not=I\se[1,k]$
was first observed by Neumann [N54b].
When all the $a_i$ are the identity element $e$
of $G$, McCoy [Mc57] obtained
the existence of a positive integer $n$
such that $x^n\in\bigcap_{i=1}^kG_i$ for all $x\in G=\bigcap_{i\in
\em}G_i$, and S. Wigzell [W95] deduced the inequality
$[G:\bigcap^k_{i=1}G_i]\ls k!\sum^k_{l=1}(-1)^{l-1}/{l!}$.
In contrast with McCoy's result, J. Backelin [B94] proved
that if a (not necessarily unitary) ring $R$ is irredundantly covered
by finitely many ideals $I_1,\ldots,I_k$ then
$$R^{k-1}:=\bg\{\sum_{j=1}^lr_{j,1}\cdots r_{j,k-1}:\,
l\in\Z^+\ \t{and}\ r_{j,1},\ldots,r_{j,k-1}\in R\bg\}$$
is contained in $\bigcap_{i=1}^kI_i$.
\medskip

 The following example shows that those inequalities in Theorem 2.1(i) are essentially sharp.
\smallskip

{\it Example} 2.1. Let $G$ be the symmetric group $S_k$ on $\{1,\ldots,k\}\ (k>1)$
and $H$ be the stabilizer of 1. Then
$$\{G_1,\ (12)G_2,\ \ldots,\ (1k)G_k\}
=\{H,\ H(12),\ H(13),\ \ldots,\ H(1k)\}$$
forms a partition of $G$,
where $G_i$ is the stabilizer of $i$ for each $i=1,\ldots,k$.
Let $\em\ne I\se[1,k]$ and $\bar I=[1,k]\sm I$.
Clearly
$\bigcap_{j\in\bar I}G_j\cong S_{|I|}$ and hence
$\l|\bigcap_{j\in\bar I}G_j\r|=|I|!$,
also $x^{\prod^{|I|}_{s=1}s}\in\bigcap^k_{i=1}G_i=\{e\}$
for all $x\in\bigcap_{j\in\bar I}G_j$.
Let $a_1=e$, and $a_i=(1i)$ for $i\in[2,k]$.
If $\sigma\in G=S_k$ and $\sigma^{-1}(1)\in I$, then
$\sigma\bigcap_{j\in\bar I}G_j\se\bigcup_{i\in I}a_iG_i$, because
for any $\tau\in\bigcap_{s\in\bar I}G_s$ and $j\in\bar I$ we have
$\sigma\tau(j)=\sigma(j)\not=1$ and hence $\sigma\tau\not\in a_jG_j$.
Note also that
 $$\sum_{i\in I}\f1{[G:G_i\cap\bigcap_{j\in\bar I}G_j]}
=\sum_{i\in I}\f{(|I|-1)!}{k!}=\f1{[G:\bigcap_{j\in\bar I}G_j]}$$
and $$\bg|\bigcup_{i\in I}G_i\bg|=\sum^{|I|}_{l=1}(-1)^{l-1}
\sum\Sb J\se I\\|J|=l\endSb\bg|\bigcap_{j\in J}G_j\bg|
=\sum^{|I|}_{l=1}(-1)^{l-1}\bi{|I|}l(k-l)!.$$

\proclaim{Lemma 2.1} Let $H$ and $K$ be subgroups of a group $G$.
If $[G:H]$ and $[G:K]$ are finite and relatively prime, then $HK$
coincides with $G$.
\endproclaim
\Proof. Use $[m,n]$ to denote the least common multiple of $m,n\in\Z^+$.
Then
$$[[G:H],[G:K]]\mid[G:H\cap K]\ \ \t{and}\ \
[G:H\cap K]\ls [G:H][G:K].$$
As $([G:H],[G:K])=1$, we have $[G:H\cap K]=[G:H][G:K]$, i.e.,
$[K:H\cap K]=[G:H]$. Since $HK$ contains exactly $[K:H\cap K]$
right cosets of $H$, by the above we must
have $HK=G$. \qed

\Remark\ 2.2. Combining Lemma 2.1 with [S01, Lemma 2.1(i)]
we find that for two subgroups $H$ and $K$ of $G$ with finite
index if $aH\cap bK=\em$ for some $a,b\in G$ then
$([G:H],[G:K])>1$. This confirms Conjecture 1.2 in the case $k=2$.
\medskip

Recall that for a subgroup $H$ of a group $G$ we use
$G/H$ to denote the set of all left cosets of $H$ in $G$ even if
$H$ may not be normal in $G$.

\proclaim{Theorem 2.2} Let $G_1,\ldots,G_k,H$ be
subgroups of a group $G$ such that either $G_1,\ldots,G_k$ are normal and $H$
is maximal in $G$, or $G_1,\ldots,G_k$ are subnormal and $H$ is
maximal normal in $G$. Suppose that $(1.1)$ is a regular system
for some $a_1,\ldots,a_k\in G$ and $[G:G_i]<\infty$ for all $i=1,\ldots,k$.
If for any $i,j=1,\ldots,k$,
 either $[G:G_i]$ is relatively prime to
$[G_i:G_i\cap G_j]$, or $G_i$ and $ G_j$ are normal in $G$ with
$G/(G_i\cap G_j)$ cyclic,
then we have
$$\{C\in G/H: C\sp a_iG_i\ \t{for some}\ i=1,\ldots,k\}
=\em\ \t{or}\  G/H\tag2.5$$
provided that $(1.1)$ is a cover of $G$ or
not all the $G_i$ are contained in $H$.
\endproclaim

\Proof. Let $I=\{1\ls  i\ls  k\colon  G_i\not\se H\}$.
If $I=[1,k]$, then the left-hand side of (2.5) is
empty. (Note that if $aH\sp  a_iG_i$ for some $a\in G$
then $a_iH=aH\sp  a_iG_i$
and hence $H\sp  G_i$.)
When $I$ is empty and (1.1) is a cover of $G$,
the left-hand side of (2.5) coincides with $G/H$
because $xH\supseteq a_iG_i$ if $x\in a_iG_i$.

Now assume that $\em\not=I\not=[1,k]$ and fix a $j\in[1,k]\sm I$.
Since $H\supseteq G_j$, $[G:H]\ls[G:G_j]<\infty$.
By Theorem 2.1(ii),
the left-hand side of (2.5) contains at least
$[(\bigcap_{i\in  I} G_i)H:H]$
left cosets of $H$. So it suffices to show that $(\bigcap_{i\in  I} G_i)H=G$.

 Let $i\in I$. As $G_i\not\se H$, we have $G_iH=HG_i=G$
 by [S01, Lemma 2.1(ii)].
Observe that $[G:H]=[G_iH:H]=[G_i:G_i\cap H]$ divides
$[G_i:G_i\cap G_j]$ and that
 $$[G:G_i\cap H]=[G:G_i][G_i:G_i\cap  H]=[G:G_i][G_iH:H]=[G:G_i][G:H].$$
If $([G:G_i],[G_i:G_i\cap G_j])=1$, then
$([G:G_i],[G:H])=1$. When $G_i,G_j$ are normal in $G$ and
$G/(G_i\cap G_j)$ is cyclic, we have
$$\align&[G/(G_i\cap G_j):G_i/(G_i\cap G_j)\cap H/(G_i\cap G_j)]
\\=&[[G/(G_i\cap G_j):G_i/(G_i\cap G_j)],[G/(G_i\cap G_j):H/(G_i\cap G_j)]],
\endalign$$
hence $[G:G_i][G:H]=[G:G_i\cap H]=[[G:G_i],[G:H]]$ and therefore
$([G:G_i],[G:H])=1$.

In light of [S01, Lemma 3.1(ii)], $[G:\bigcap_{i\in I} G_i]$ divides
$\prod_{ i\in  I}[G:G_i]$. So
$([G:\bigcap_{i\in  I} G_i],[G:H])=1$ by the above. Hence
$\l(\bigcap_{i\in I}G_i\r)H=G$ by Lemma 2.1,
and this concludes the proof. \qed

\proclaim{Corollary 2.1} Let $(1.1)$ be a regular cover of a group $G$
such that for any $i,j=1,\ldots,k$
either $G_i$ and $G_j$ are normal in $G$ with $G/(G_i\cap G_j)$ cyclic,
or $G_i$ and $G_j$ are subnormal in $G$ with
$([G:G_i],[G_i:G_i\cap G_j])=1$.
If $G_j\not=G$, then for some $i\not=j$ with $G_i\not=G$ we have
$a_iG_i\cap a_jG_j=\em$.
\endproclaim
\Proof. Suppose that $G_j\ne G$ where $1\ls j\ls k$.
As $G_j$ is subnormal and of finite index in $G$
it is contained in a proper maximal normal subgroup $H$ of $G$.
Clearly $a_jG_j\se a_jH\subset G$. Choose $x\in G\sm a_jH$.
By Theorem 2.2, $xH\sp a_iG_i$ for some $i\in[1,k]$.
Since $xH$ is disjoint from $a_jH$, we have $a_iG_i\cap a_jG_j=\em$.
Obviously $i\ne j$ and $G_i\ne G$. We are done. \qed

\heading{3. Proofs of Theorems 1.3-1.5}\endheading

\noindent{\it Proof of Theorem 1.3}.
We use induction on the finite index $[G:\bigcap^k_{i=1}G_i]$.

If $[G:\bigcap^k_{i=1}G_i]=1$, then $G_1=\cs=G_k=G$, $m(\Cal A)=k$
and $d(G,\bigcap^k_{i=1}G_i)=0$, so (1.9)
holds trivially.

Now let's proceed to the induction step and assume that
$[G:\bigcap^k_{i=1}G_i]>1$.

The following observations are important in our induction step.
Let $K$ be any subgroup of $G$. Then $G_i\cap K$ is subnormal in $K$
for any $i=1,\ldots,k$.
Let $i,j\in[1,k]$. If $([G:G_i],[G_i:G_i\cap G_j])=1$,
then $[K:G_i\cap K]$ is relatively prime to $[G_i\cap K:G_i\cap G_j\cap K]$
because $[K:G_i\cap K]\mid [G:G_i]$
and $[G_i\cap K:(G_i\cap G_j)\cap (G_i\cap K)]\mid [G_i:G_i\cap G_j]$
by [S01, Lemma 3.1]. If both $G_i$ and $G_j$ are normal in $G$
with $G/(G_i\cap G_j)$ cyclic, then both $G_i\cap K$ and $G_j\cap K$
are normal in $K$
and $K/(G_i\cap G_j\cap K)\cong (G_i\cap G_j)K/(G_i\cap G_j)$ is cyclic.

Suppose $G_{i_0}\ne G$ where $1\ls i_0\ls k$.
As $G_{i_0}$ is subnormal in $G$ there is a maximal normal subgroup
$H$ of $G$ with $G_{i_0}\se H\subset G$.
Let $I_0$ be a minimal subset of $[1,k]$
such that $|\{i\in I_0\colon x\in a_iG_i\}|\gs m(\Cal A)$ for all $x\in H$.
Clearly such an $I_0$ exists and $a_iG_i\cap H\ne\em$ for all $i\in I_0$.
Observe that the nonempty system $\Cal A_0=\{a_iG_i\cap H\}_{i\in I_0}$
forms a minimal $m(\Cal A)$-cover
of $H$ by left cosets of subnormal subgroups $G_i\cap H$ ($i\in I_0$) of $H$.
Also,
$$\[H:\bigcap_{i\in I_0}(G_i\cap H)\]\ls\[H:\bigcap^k_{i=1}G_i\cap H\]
=\[H:\bigcap^k_{i=1}G_i\]<\[G:\bigcap^k_{i=1}G_i\].$$
Since a minimal $m(\Cal A)$-cover is also a regular cover,
by the induction hypothesis we have
$$|I_0|\gs m(\Cal A_0)+d\(H,\bigcap_{i\in I_0}(G_i\cap H)\)
= m(\Cal A)+d\(H,H\cap\bigcap_{i\in I_0}G_i\).$$

In the case  $I_0\not=[1,k]$, by Theorem 2.1(i) there exists a $g_1\in G$
such that
$$\bigcup^k\Sb j=1\\j\not\in I_0\endSb a_jG_j\sp g_1\bigcap_{i\in I_0}G_i,$$
hence we can choose a minimal $I_1\se[1,k]$ with $I_1\cap I_0=\em$ such that
$$\bigcup_{j\in I_1}a_jG_j
\cap g_1H\sp g_1\l(\bigcap_{i\in I_0}G_i\cap H\r).$$
Since $\Cal A_1=\{g_1^{-1}a_jG_j\cap\bigcap_{i\in I_0}G_i\cap H\}_{j\in I_1}$
forms a minimal cover of $\bigcap_{i\in I_0}G_i\cap H$
by left cosets of subnormal subgroups
$G_j\cap\bigcap_{i\in I_0}G_i\cap H$ ($j\in I_1)$ of
$\bigcap_{i\in I_0}G_i\cap H$, again by the induction
hypothesis we have
$$\aligned|I_1|\gs& m(\Cal A_1)
+d\(\bigcap_{i\in I_0}G_i\cap H,\bigcap_{j\in I_1}
\(G_j\cap\bigcap_{i\in I_0}G_i\cap H\)\)
\\=&1+d\(H\cap\bigcap_{i\in I_0}G_i,H\cap
\bigcap_{i\in I_0\cup I_1}G_i\).
\endaligned$$

Continue the above procedure until we obtain a partition
 $\{I_s\}^n_{s=0}$ (with $0\ls n< k$)
of $[1,k]$ such that
for each $s\in[1,n]$ there is a $g_s\in G$
with $a_iG_i\cap g_sH\ne\em$ for all $i\in I_s$
and
$$|I_s|\gs1+d\(H\cap\bigcap_{i\in I_0\cup\cs\cup I_{s-1}}G_i,
H\cap\bigcap_{i\in I_0\cup\cs\cup I_s}G_i\).$$
Observe that
$$\aligned k=\sum^n_{s=0}|I_s|\gs&m(\Cal A)+d\(H,H\cap\bigcap_{i\in I_0}G_i\)
\\&\ +\sum_{0<s\ls n}\(1+d\(H\cap\bigcap_{i\in I_0\cup\cs\cup I_{s-1}}G_i,
H\cap\bigcap_{i\in I_0\cup\cs\cup I_s}G_i\)\)
\\=&m(\Cal A)+n+d\(H,H\cap\bigcap_{i\in I_0\cup\cs\cup I_n}G_i\)
\\=&m(\Cal A)+n+d\(H,\bigcap^k_{i=1}G_i\).\endaligned$$

 As $a_{i_0}H\sp a_{i_0}G_{i_0}$, by Theorem 2.2 each $C\in G/H$
 contains $a_iG_i$
for some $i\in[1,k]=I_0\cup\cs\cup I_n$.
Recall that $a_iG_i\cap g_0H\not=\em$ for all $i\in I_0$ where $g_0=e$. For any $g\in G$, if
$gH\sp a_iG_i$ and $i\in I_s$ then $gH\cap g_sH\sp a_iG_i\cap g_sH\not=\em$
and hence $gH=g_sH$.
Therefore
$$|G/H|=|\{g_sH\colon 0\ls s\ls n\}|\ls n+1$$ and hence
$$\aligned k-m(\Cal A)\gs&n+d\(H,\bigcap^k_{i=1}G_i\)\gs |G/H|-1
+d\(H,\bigcap^k_{i=1}G_i\)
\\=&d(G,H)+d\(H,\bigcap^k_{i=1}G_i\)=d\(G,\bigcap^k_{i=1}G_i\).\endaligned$$
This completes the induction proof. \qed

\Remark\ 3.1. In view of Theorem 2.2, by modifying the
complicated proof of [S90, Theorem 8],
we can obtain the following result for minimal $m$-covers
({\it not} for regular covers):
Let (1.1) be a minimal $m$-cover of a group $G$, and suppose that
for any $i,j=1,\ldots,k$ either $G_i$ and $G_j$
are subnormal in $G$ with $([G:G_i],[G_i:G_i\cap G_j])=1$,
or $G_i$ and $G_j$ are normal in $G$ with $G/(G_i\cap G_j)$
cyclic. Then, for each subgroup $K$ of $G$ with
$I(K)=\{1\ls i\ls k:\,K\not\se G_i\}\not=\em$, we can find
an $r\in I(K)$ and $x_i\in K\sm G_i$ (for $i\in I(K)\sm\{r\}$)
such that
$$\bg|\bg\{x_i\(K\cap\bigcap_{s=1}^kG_s\):\, i\in I(K)\sm\{r\}\bg\}\bg|
\gs d\(K,K\cap\bigcap_{s=1}^kG_s\).$$
If $K\supseteq G_j$ for some $j\in[1,k]$,
or $K\supseteq H=(\bigcap_{s=1}^kG_s)_G$
and $K/H$ is a Hall subgroup of $G/H$, then
$I(K)=\{1\ls i\ls k:\,[G:G_i]\nmid[G:K]\}$ by Lemma 3.1 given later.

\medskip
\noindent{\it Proof of Theorem 1.4}.
Let $M$ be any perfect normal subgroup of $G$.
Fix $x\in G$ and
set $M_x=\langle M,x\rangle=\langle x\rangle M$.
Since $M_x/M\cong\langle x\rangle/(\langle x\rangle\cap M)$ is cyclic
and hence abelian, $M_x'
\se M=M'\se M'_x$ and therefore $M'_x=M$.

As $\{G_i\}_{i=1}^k$ is an $m$-cover of $G$, we can choose a
minimal $I_x\se[1,k]$ such that $|I_x^*(g)|\gs m$ for all $g\in
M_x$. Clearly $\Cal A_x=\{G_i\cap M_x\}_{i\in I_x}$ forms a
minimal $m$-cover of $M_x$ by subnormal subgroups of $M_x$.
We claim that $J=\{i\in I_x:\,G_i\not\supseteq M_x\}$ is empty.
Assume on the contrary
that $J\not=\em$. Then $|\{i\in I_x:\, G_i\supseteq M_x\}|<m$.
Since $\{G_i\cap M_x\}_{i\in J}$ is a cover of
$M_x$ by proper subnormal subgroups,  $M_x$ possesses a
cover by proper normal subgroups. Applying a result of Brodie,
Chamberlain and Kappe [BCK], we find that $M_x$ has a normal
subgroup $H$ such that $M_x/H\cong C_p\times C_p$ for some prime
$p$. As $M_x/H$ is abelian, we have $M=M_x'\se H$ and hence
$M_x/H\cong (M_x/M)/(H/M)$ is cyclic. This contradiction shows
that our claim is true.

Let $I=\bigcup_{x\in G}I_x$.
For each $x\in G$, clearly $|I^*(x)|\gs|I_x^*(x)|\gs m$.
Since $\{G_i\}_{i=1}^k$ is a minimal $m$-cover of $G$,
we must have $I=[1,k]$. For any $i\in[1,k]$,
there is an $x\in G$ such that $i\in I_x$ and hence
$G_i\supseteq M_x\supseteq M$.

By the above, $F=\bigcap_{i=1}^kG_i$ contains
any perfect normal subgroup of $G$.
As $\{G_i/F_G\}_{i=1}^k$ is a minimal $m$-cover of
$G/F_G$, $F/F_G=\bigcap_{i=1}^kG_i/F_G$ contains
any perfect normal subgroup of the finite group $G/F_G$.
In light of parts (i)and (v) of Theorem 3.1 of [S01],
there exists a composition series from $F/F_G$ to $G/F_G$ whose factors
are of prime orders.
Thus there is also a composition series from $F$ to $G$ whose factors
are of prime orders, and hence
$F$ contains all perfect subgroups of $G$ by [S01, Theorem 3.1].
This concludes our proof. \qed

\proclaim{Lemma 3.1}
Let $H$ and $K$ be subgroups of a group
$G$ with finite index.
Suppose that $H$ or $K$ is subnormal in $G$, and $([G:K],[K:H\cap K])=1$. Then
$$[G:H\cap K]=[[G:H],[G:K]];\tag3.1$$
hence
$$H\supseteq K\iff[G:H]\mid[G:K],\quad\t{and}\quad
K\supseteq H\iff[G:K]\mid[G:H].$$
\endproclaim
\Proof. As $[G:K]$ is relatively prime to $[K:H\cap K]$,
$$[G:\hk]=[G:K][K:\hk]=[[G:K],[K:\hk]].$$
By [S01, Lemma 3.1(i)], $[K:\hk]\mid[G:H]$ and therefore $[G:\hk]\mid[\h,[G:K]]$.
On the other hand, $[G:\hk]$
is obviously divisible by both $\h$ and $[G:K]$.
So we have (3.1).

Observe that
$$\aligned H\supseteq K&\iff [G:\hk]/[G:K]=[K:\hk]=1
\\&\iff [\h,[G:K]]=[G:K]\ \t{(i.e.,}\ [G:H]\mid[G:K]).\endaligned$$
Similarly, $K\supseteq H$ if and only if $[G:K]\mid[G:H]$.
This ends the proof. \qed

\medskip
\noindent{\it Proof of Theorem 1.5}.
For any $\em\not=J\se[1,k]$, we assert that
$[G:\bigcap_{j\in J}G_j]$ is the least common multiple $[n_j]_{j\in J}$
of those $n_j=[G:G_j]$ with $j\in J$.
When $|J|=1$, this is trivial.
Now assume $|J|>1$ and
$[G:\bigcap_{j\in J\sm\{j_0\}}G_j]=[n_j]_{j\in J\sm\{j_0\}}$,
where $j_0=\max J$. Observe that
$[G_{j_0}:G_{j_0}\cap\bigcap_{j\in J\sm\{j_0\}}G_j]$ divides $[G_{j_0}:H]$ and
hence it is relatively prime to $[G:G_{j_0}]$. With the help of Lemma 3.1, we find that
$$\align\[G:\bigcap_{j\in J}G_j\]=&\[[G:G_{j_0}],\[G:\bigcap_{j\in J\sm\{j_0\}}G_j\]\]
\\=&[n_{j_0},[n_j]_{j\in J\sm\{j_0\}}]=[n_j]_{j\in J}.
\endalign$$
This proves the assertion by induction.

In view of the above and the inclusion-exclusion principle,
$$\align\[\bigcup_{i=1}^kG_i:H\]
=&\sum_{\em\not=J\se[1,k]}(-1)^{|J|-1}\[\bigcap_{j\in J}G_j:H\]
\\=&\sum_{\em\not=J\se[1,k]}(-1)^{|J|-1}\f{[G:H]}{[n_j]_{j\in J}}
\\=&\sum_{\em\not=J\se[1,k]}(-1)^{|J|-1}
\bg|\bg\{0\ls x<[G:H]:\,x\in\bigcap_{j\in J}n_j\Z\bg\}\bg|
\\=&\sum_{\em\not=J\se[1,k]}(-1)^{|J|-1}\bg|\bigcap_{j\in J}X_j\bg|
=\bg|\bigcup_{i=1}^kX_i\bg|
\endalign$$
where $X_i=\{0\ls x<[G:H]:\,x\in n_i\Z\}$ for $i=1,\ldots,k$.
Applying [S04b, Theorem 3.1], we finally obtain that
$$\align\[\bigcup_{i=1}^ka_iG_i:H\]\gs&\bg|\bg\{0\ls x<[G:H]:
\,x\in\bigcup_{i=1}^k n_i\Z\bg\}\bg|
\\=&\bg|\bigcup_{i=1}^kX_i\bg|=\[\bigcup_{i=1}^kG_i:H\].
\endalign$$

 The proof of Theorem 1.5 is now complete. \qed

 \Remark\ 3.2. If $a_1G_1,\ldots,a_kG_k$ are left cosets
 in a finite cyclic group $G$, then we also have
 (1.12) because $[G:H\cap K]=[[G:H],[G:K]]$
 (i.e., $|H\cap K|=(|H|,|K|)$) for any
 subgroups $H$ and $K$ of $G$.

\Ack.  The paper was finished during the author's stay
at the Institute of Camille Jordan (Univ. Lyon-I, France)
as a visiting professor, so he would like to
thank Prof. J. Zeng for the invitation and
hospitality.

\enddocument